\magnification\magstep1
%\advance\vsize2truecm
\baselineskip15pt
%\newread\AUX\immediate\openin\AUX=\jobname.aux
\def\ref#1{\expandafter\edef\csname#1\endcsname}
%
%\ifeof\AUX\immediate\write16{\jobname.aux gibt es nicht!}\else
%\input \jobname.aux
%\fi\immediate\closein\AUX
%
\ref {Introduction}{1}
\ref {Definitions}{2}
\ref {AAA}{Theorem\penalty 10000\ 2.1}
\ref {BBB}{Theorem\penalty 10000\ 2.2}
\ref {Hecke Algebra}{3}
\ref {Self}{Lemma\penalty 10000\ 3.1}
\ref {Recursion}{4}
\ref {Haction}{Theorem\penalty 10000\ 4.2}
\ref {Recurs}{Corollary\penalty 10000\ 4.4}
\ref {Integrality}{5}
\ref {Recurs2}{Theorem\penalty 10000\ 5.1}
\ref {Weak Integrality}{Corollary\penalty 10000\ 5.2}
\ref {Main1}{Theorem\penalty 10000\ 5.4}
\ref {OpInt}{Lemma\penalty 10000\ 5.5}
\ref {Phiop}{Lemma\penalty 10000\ 5.6}
\ref {q=0}{Corollary\penalty 10000\ 5.7}
\ref {Symmetric}{6}
\ref {MainS}{Theorem\penalty 10000\ 6.1}
\ref {Jack}{7}
\ref {References}{8}
\ref {C1}{[C1]}
\ref {LS}{[LS]}
\ref {M3}{[M1]}
\ref {M2}{[M2]}
\ref {M1}{[M3]}
\ref {Op}{[O1]}
\def\today{\number\day.~\ifcase\month\or
  Januar\or Februar\or M\"arz\or April\or Mai\or Juni\or
  Juli\or August\or September\or Oktober\or November\or Dezember\fi
  \space\number\year}
\font\sevenex=cmex7
%\font\sevenex=cmex10 scaled 700
\scriptfont3=\sevenex
%\font\fiveex=cmex7 scaled 714
\font\fiveex=cmex10 scaled 500
\scriptscriptfont3=\fiveex

\def\phi{\varphi}
\def\epsilon{\varepsilon}
\def\theta{\vartheta}

%\font\lams=lams3
%\def\uinto{\lower1.7pt\hbox{%
%\vbox{\offinterlineskip
%\hbox{\lams\char"7A}%
%\hbox{\vbox to 7.5pt{\leaders\vrule width0.2pt\vfill}%
%\kern-.3pt\hbox{\lams\char"76}}}}}
\def\uauf{\lower1.7pt\hbox to 3pt{%
\vbox{\offinterlineskip
\hbox{\vbox to 8.5pt{\leaders\vrule width0.2pt\vfill}%
\kern-.3pt\hbox{\lams\char"76}\kern-0.3pt%
$\raise1pt\hbox{\lams\char"76}$}}\hfil}}
%%%%%%%%%%%%%%%%%%
% Makros f"ur Querverweise:
\def\cite#1{\expandafter\ifx\csname#1\endcsname\relax
{\bf?}\immediate\write16{#1 ist nicht definiert!}\else\csname#1\endcsname\fi}
\def\expandwrite#1#2{\edef\next{\write#1{#2}}\next}
\def\neverexpand{\noexpand\noexpand\noexpand}
\def\strip#1\ {}
\def\ncite#1{\expandafter\ifx\csname#1\endcsname\relax
{\bf?}\immediate\write16{#1 ist nicht definiert!}\else
\expandafter\expandafter\expandafter\strip\csname#1\endcsname\fi}
\newwrite\AUX
\immediate\openout\AUX=\jobname.aux
%%%%%%%%%%%%%%%%%%%%%
\newcount\Abschnitt\Abschnitt0
\def\beginsection#1. #2 \par{\advance\Abschnitt1%
\vskip0pt plus.10\vsize\penalty-250
\vskip0pt plus-.10\vsize\bigskip\vskip\parskip
\edef\TEST{\number\Abschnitt}
\expandafter\ifx\csname#1\endcsname\TEST\relax\else
\immediate\write16{#1 hat sich geaendert!}\fi
\expandwrite\AUX{\neverexpand\ref{#1}{\TEST}}
\leftline{\bf\number\Abschnitt. \ignorespaces#2}%
\nobreak\smallskip\noindent\SATZ1}
%%%%%%%%%%%%%%%%%%
\def\Proof:{\par\noindent{\it Proof:}}
\def\Remark:{\ifdim\lastskip<\medskipamount\removelastskip\medskip\fi
\noindent{\bf Remark:}}
\def\Remarks:{\ifdim\lastskip<\medskipamount\removelastskip\medskip\fi
\noindent{\bf Remarks:}}
\def\Definition:{\ifdim\lastskip<\medskipamount\removelastskip\medskip\fi
\noindent{\bf Definition:}}
\def\Example:{\ifdim\lastskip<\medskipamount\removelastskip\medskip\fi
\noindent{\bf Example:}}
%%%%%%%%%%%%%%%%
\newcount\SATZ\SATZ1
\def\proclaim #1. #2\par{\ifdim\lastskip<\medskipamount\removelastskip
\medskip\fi
\noindent{\bf#1.\ }{\it#2}\Par
\ifdim\lastskip<\medskipamount\removelastskip\goodbreak\medskip\fi}
\def\Aussage#1{%
\expandafter\def\csname#1\endcsname##1.{\ifx?##1?\relax\else
\edef\TEST{#1\penalty10000\ \number\Abschnitt.\number\SATZ}
\expandafter\ifx\csname##1\endcsname\TEST\relax\else
\immediate\write16{##1 hat sich geaendert!}\fi
\expandwrite\AUX{\neverexpand\ref{##1}{\TEST}}\fi
\proclaim {\number\Abschnitt.\number\SATZ. #1\global\advance\SATZ1}.}}
\Aussage{Theorem}
\Aussage{Proposition}
\Aussage{Corollary}
\Aussage{Lemma}
%%%%%%%%%%%%%%%%
\font\la=lasy10
\def\strich{\hbox{$\vcenter{\hbox
to 1pt{\leaders\hrule height -0,2pt depth 0,6pt\hfil}}$}}
\def\dashedrightarrow{\hbox{%
\hbox to 0,5cm{\leaders\hbox to 2pt{\hfil\strich\hfil}\hfil}%
\kern-2pt\hbox{\la\char\string"29}}}

\def\Bindestrich{\penalty10000-\hskip0pt}
\let\_=\Bindestrich
\def\.{{\sfcode`.=1000.}}
%%%%%%%%%%%%%%%%%%%%%%%%%%%%%%%%%%%%

\def\Par{\par}
\def\:={\mathrel{\raise0,9pt\hbox{.}\kern-2,77779pt
\raise3pt\hbox{.}\kern-2,5pt=}}
\def\=:{\mathrel{=\kern-2,5pt\raise0,9pt\hbox{.}\kern-2,77779pt
\raise3pt\hbox{.}}} 

\def\pfeil{\rightarrow}

\def\Ugleich{\hbox{$\cup$\kern.5pt\vrule depth -0.5pt}}
\def\|#1|{\mathop{\rm#1}\nolimits}
\def\<{\langle}
\def\>{\rangle}
\let\Times=\times
\def\times{\mathop{\Times}}
\let\Otimes=\otimes
\def\otimes{\mathop{\Otimes}}
%%%%%%%%%%%%%%%%%%%%%%%%%%%%%%%%%
%Laden von Fonts:
\catcode`\@=11
\def\hex#1{\ifcase#1 0\or1\or2\or3\or4\or5\or6\or7\or8\or9\or A\or B\or
C\or D\or E\or F\else\message{Warnung: Setze hex#1=0}0\fi}
\def\fontdef#1:#2,#3,#4.{%
\alloc@8\fam\chardef\sixt@@n\FAM
\ifx!#2!\else\expandafter\font\csname text#1\endcsname=#2
\textfont\the\FAM=\csname text#1\endcsname\fi
\ifx!#3!\else\expandafter\font\csname script#1\endcsname=#3
\scriptfont\the\FAM=\csname script#1\endcsname\fi
\ifx!#4!\else\expandafter\font\csname scriptscript#1\endcsname=#4
\scriptscriptfont\the\FAM=\csname scriptscript#1\endcsname\fi
\expandafter\edef\csname #1\endcsname{\fam\the\FAM\csname text#1\endcsname}
\expandafter\edef\csname hex#1fam\endcsname{\hex\FAM}}
\catcode`\@=12 

%%%%%%%%%%%%%%%%%%%%%%%%%%%%%%%%%
\fontdef Ss:cmss10,,.
\fontdef Fr:eufm10,eufm7,eufm5.

			%Hier aufpassen!!!

\def\fm{{\Fr m}}

\def\fQ{{\Fr Q}}

\newread\AUXX
\immediate\openin\AUXX=msxym.tex
\ifeof\AUXX
\fontdef bbb:msbm10,msbm7,msbm5.
\fontdef mbf:cmmib10,cmmib7,.
\else
\fontdef bbb:msym10,msym7,msym5.
\fontdef mbf:cmmib10,cmmib10 scaled 700,.
\fi
\immediate\closein\AUXX
\def\CC{{\bbb C}}

\def\NN{{\bbb N}}
\def\QQ{{\bbb Q}}

\def\ZZ{{\bbb Z}}

\def\cE{{\cal E}}
\def\cJ{{\cal J}}
\def\cP{{\cal P}}

\mathchardef\leer=\string"0\hexbbbfam3F
\mathchardef\subsetneq=\string"3\hexbbbfam24
\mathchardef\semidir=\string"2\hexbbbfam6E
\mathchardef\dirsemi=\string"2\hexbbbfam6F
\let\OL=\overline
\def\overline#1{{\hskip1pt\OL{\hskip-1pt#1\hskip-1pt}\hskip1pt}}
\def\Aq{{\overline{A}}}

\def\fQ{{\overline{f}}}%<--                    Aufpassen  

\def\gq{{\overline{g}}}
\def\Hq{{\overline{H}}}

\def\qq{{\overline{q}}}

\def\tq{{\overline{t}}}

\def\xq{{\overline{x}}}

%
%%%%%%%%%%%%
% Displayroutine
\abovedisplayskip 9.0pt plus 3.0pt minus 3.0pt
\belowdisplayskip 9.0pt plus 3.0pt minus 3.0pt
\newdimen\Grenze\Grenze2\parindent\advance\Grenze1em
\newdimen\Breite
\newbox\DpBox
\def\NewDisplay#1$${\Breite\hsize\advance\Breite-\hangindent
\setbox\DpBox=\hbox{\hskip2\parindent$\displaystyle{#1}$}%
\ifnum\predisplaysize<\Grenze\abovedisplayskip\abovedisplayshortskip
\belowdisplayskip\belowdisplayshortskip\fi
\global\futurelet\nexttok\WEITER}
\def\WEITER{\ifx\nexttok\qed\expandafter\leftQEDdisplay
\else\leftdisplay\fi}
\def\leftdisplay{\hskip-\hangindent\leftline{\box\DpBox}$$}
\def\leftQEDdisplay{\hskip-\hangindent
\line{\copy\DpBox\hfill\lower\dp\DpBox\copy\QEDbox}%
\belowdisplayskip0pt$$\bigskip\let\nexttok=}
\everydisplay{\NewDisplay}
%%%%%%%%%%%%
\newbox\QEDbox
\newbox\nichts\setbox\nichts=\vbox{}\wd\nichts=2mm\ht\nichts=2mm
\setbox\QEDbox=\hbox{\vrule\vbox{\hrule\copy\nichts\hrule}\vrule}
\def\qed{\leavevmode\unskip\hfil\null\nobreak\hfill\copy\QEDbox\medbreak}
%%%%%%%%%%%%%%
\newdimen\HIindent
\newbox\HIbox
\def\setHI#1{\setbox\HIbox=\hbox{#1}\HIindent=\wd\HIbox}
\def\HI#1{\par\hangindent\HIindent\hangafter=0\noindent\leavevmode
\llap{\hbox to\HIindent{#1\hfil}}\ignorespaces}
%%%%%%%%%%%%%%
\def\rho{\varrho}

\def\lamq{{\overline\lambda}}
\fontdef Ss:cmss10,,.
\font\BF=cmbx10 scaled \magstep2
\font\CSC=cmcsc10 %scaled \magstephalf
\baselineskip15pt
{\baselineskip1.5\baselineskip\rightskip0pt plus 5truecm
\leavevmode\vskip0truecm\noindent
\BF Integrality of Two Variable Kostka Functions

}
\vskip1truecm
\leftline{{\CSC Friedrich Knop}%
\footnote*{\rm Partially supported by a grant of the NSF}}
\leftline{Department of Mathematics, Rutgers University, New Brunswick NJ
08903, USA}
\leftline{knop@math.rutgers.edu}
\vskip1truecm
\beginsection Introduction. Introduction

Macdonald defined in \cite{M3} a remarkable class of symmetric
polynomials $P_\lambda(x;q,t)$ which depend on two parameters and
interpolate between many families of classical symmetric polynomials. For
example $P_\lambda(x;t)=P_\lambda(x;0,t)$ are the Hall\_Littlewood
polynomials which themselves specialize for $t=0$ to Schur functions
$s_\lambda$. Also Jack polynomials arise by taking $q=t^\alpha$ and
letting $t$ tend to $1$.

The Hall\_Littlewood polynomials are orthogonal with
respect to a certain scalar product $\<\cdot,\cdot\>_t$. The scalar
products $K_{\lambda\mu}=\<s_\lambda,s_\mu\>_0$ are known as Kostka
numbers. Since it has an interpretation as a weight multiplicity of a
$GL_n$\_representation, it is a natural number. Also the scalar products 
$K_{\lambda\mu}(t)=\<P_\lambda(x;t),s_\mu\>_t$ are important
polynomials in $t$. Their coefficients have been proven to be natural
numbers by Lascoux\_Sch\"utzenberger \cite{LS}.

Macdonald conjectured in \cite{M3} that even the scalar
products 
$$
K_{\lambda\mu}(q,t)=\<P_\lambda(x;q,t),s_\mu\>_t
$$
are polynomials in $q$ and $t$ with non\_negative integers as
coefficients. In this note we prove that $K_{\lambda\mu}(q,t)$ is at
least a polynomial. The positivity seems to be much deeper and is, to my
knowledge, unsolved up to know.

Our main tool has also been introduced by Macdonald. Following a
construction of Opdam \cite{Op} in the Jack polynomial case, he
constructed a family $E_\lambda$ of non\_symmetric polynomials. Here
the indexing set is now all of $\NN^n$. They have properties which are
very similar to those of the symmetric functions, but in a sense they are
easier to work with. In particular, we exhibit a very simple recursion
formula (a ``creation operator'') in terms of Hecke operators for them.
This formula enables us to prove an analogous conjecture for
these non\_symmetric polynomials, at least what polynomiality concerns.
At the end, we obtain our main result by symmetrization.

The proof follows very closely a proof of an analogous conjecture for
Jack polynomials. In this case we were even able to settle positivity.
This will appear along with a new combinatorial formula for Jack
polynomials as a joint paper with S. Sahi. The main difference to the
Jack polynomial case is, of course, the appearence of Hecke operators.
Furthermore, we introduce a certain basis which is the non\_symmetric
analogue of Hall\_Littlewood polynomials. But the main point is again a
very special creation operator which in the present case is $$
\Phi f(z_1,\ldots,z_n)=z_nf(q^{-1}z_n,z_1,\ldots,z_{n{-}1}).
$$
Its existence is extremely connected to the fact that we are working
with the root system of $GL_n$ as opposed to $SL_n$. Therefore, I don't
think that the techniques presented here generalize to arbitrary other
root systems.

\beginsection Definitions. Definitions

Let $\Lambda:=\NN^n$ and $\Lambda^+$ be the subset of all partitions.
For $\lambda=(\lambda_i)\in\Lambda$ we put
$|\lambda|:=\sum_i\lambda_i$ and
$l(\lambda):=\|max|\{i\mid\lambda_i\ne0\}$ (with $l(0):=0$).
To each $\lambda\in\Lambda$ corresponds a monomial $z^\lambda$.

There is a (partial) order relation on
$\Lambda$. First, recall the usual order on the set
$\Lambda^+$: we say
$\lambda\ge\mu$ if $|\lambda|=|\mu|$ and
$$
\lambda_1+\lambda_2+\ldots+\lambda_i\ge\mu_1+\mu_2+\ldots+\mu_i
\quad\hbox{for all }i=1,\ldots,n.
$$
This order relation is extended to all of $\Lambda$ as follows.
Clearly, the symmetric group $W=S_n$ on $n$ letters acts on $\Lambda$
and for every $\lambda\in\Lambda$ there is a unique partition
$\lambda^+$ in the orbit $W\lambda$. For all permutations $w\in W$ with
$\lambda=w\lambda^+$ there is a unique one, denoted by $w_\lambda$, of
minimal length. We define $\lambda\ge\mu$ if either $\lambda^+>\mu^+$
or $\lambda^+=\mu^+$ and $w_\lambda\le w_\mu$ in the Bruhat order of
$W$. In particular, $\lambda^+$ is the unique {\it maximum\/} of
$W\lambda$.

Let $k$ be the field $\CC(q,t)$, where $q$ and $t$ are formal variables
and let $\cP:=k[z_1,\ldots,z_n]$ be the polynomial ring and
$\cP'=k[z_1,z_1^{-1},\ldots,z_n,z_n^{-1}]$ the Laurent polynomial ring
over $k$. There are involutory automorphisms $\iota:x\mapsto\xq$ of
$k/\CC$ with $\qq=q^{-1}$, $\tq:=t^{-1}$, and a $\iota$\_semilinear
involution $\cP'\mapsto\cP':f\mapsto\fQ$ with $\overline{z_i}=z_i^{-1}$.
For $f\in\cP$ let $[f]_1$ be its constant term. Fix an $r\in\NN$ and put
$t=q^r$. Then Cherednik (see \cite{C1}, \cite{M1}~\S5) defines a certain
Laurent polynomial $\delta_r(x;t)$ and a pairing on $\cP'$ which is
Hermitian with respect to $\iota$ by
$$
\<f,g\>:=[\delta_r f \gq]_1.
$$
Non\_symmetric Macdonald polynomials are defined by the following
theorem.

\Theorem AAA. {\rm(\cite{M1}~\S6)} For every $\lambda\in\Lambda$ there
is a unique polynomial $E_\lambda(z;q,t)\in\cP$ satisfying
\item{i)}$E_\lambda(z)=z^\lambda+\sum_{\mu\in\Lambda:\mu<\lambda}
c_{\lambda\mu}(q,t)z^\mu$ and
\item{ii)}$\<E_\lambda(z),z^\mu\>_r=0$ for all $\mu\in\Lambda$ with
$\mu<\lambda$ and almost all $r\in\NN$.\Par
\noindent Moreover, the collection
$\{E_\lambda\mid\lambda\in\Lambda\}$ forms a $k$\_linear basis of
$\cP$.

The symmetric group $W$ acts also on $\cP$ in the obvious way. We are
going to define a basis of $\cP^W$, the algebra of symmetric functions,
which is parametrized by $\Lambda^+$. One basis is already given by the
monomial symmetric functions $m_\lambda$, $\lambda\in\Lambda^+$. 
Next, we define the symmetric Macdonald polynomials:

\Theorem BBB. {\rm (\cite{M1}~1.5} For every $\lambda\in\Lambda^+$
there is a unique symmetric polynomial $J_\lambda(z;q,t)\in\cP^W$
satisfying
\item{i)}$J_\lambda(z)=m_\lambda+\sum_{\mu\in\Lambda^+:\mu<\lambda}
c'_{\lambda\mu}(q,t)m_\mu$ and
\item{ii)}$\<J_\lambda(z),m_\mu\>_r=0$
for all $\mu\in\Lambda^+$ with $\mu<\lambda$ and almost all
$r\in\NN$.\Par \noindent
Moreover, the collection $\{J_\lambda\mid\lambda\in\Lambda\}$ forms an
$k$\_ linear basis of $\cP^W$.

\beginsection Hecke Algebra. The Hecke algebra

The scalar product above is {\it not\/} symmetric in the variables
$z_i$. Therefore, we define operators which replace the usual
$W$\_action. Let $s_i\in W$ be the $i$\_th simple reflection. First,
we define the operators $N_i:=(z_i-z_{i{+}1})^{-1}(1-s_i)$ and then
$$
H_i:=s_i-(1-t)N_iz_i\qquad\Hq_i:=s_i-(1-t)z_{i{+}1}N_i.
$$
They satisfy the relations
$$
H_i-\Hq_i=t-1;\qquad H_i\Hq_i=t.
$$
This means that both $H_i$ and $-\Hq_i$ solve the equation
$(x+1)(x-t)=0$.
Also the braid relations hold
$$
\eqalign{
H_iH_{i{+}1}H_i&=H_{i{+}1}H_iH_{i{+}1}\qquad i=1,\ldots, n-2\cr
H_iH_j&=H_jH_i\qquad\qquad|i-j|>1\cr}
$$
This means that the algebra generated by the $H_i$ is a Hecke
algebra of type $A_{n-1}$. For all this see \cite{M1}. For
compatibility, let me remark that our parameter $t$ is $t^2$ in \cite{M1}
and our $H_i$ is $tT_i$ there. Furthermore, the simple roots (or rather
their exponentials) are $z_{i{+}1}/z_i$.

The connection with the $W$\_action is that we get the same set of
invariants in the following sense: $f\in \cP^W$ if and only if
$\Hq_i(f)=f$ for all $i$ if and only if $H_i(f)=t\,f$ for all $i$.

The braid relations imply that $H_w:=H_{i_1}\ldots H_{i_k}$ is well
defined if $w=s_{i_1}\ldots s_{i_k}\in W$ is a reduced decomposition
and similarly for $\Hq_w$. The following relations hold:
$$
z_{i{+}1}H_i=\Hq_iz_i;\quad H_iz_{i{+}1}=z_i\Hq_i\qquad
i=1,\ldots,n-1.
$$
Now, we introduce the operator $\Delta$ by
$$
\Delta f(z_1,\ldots,z_n):=f(q^{-1}z_n,z_1,\ldots,z_{n{-}1}).
$$
The following relations are easily checked
$$
\eqalign{
\Delta z_{i{+}1}&=z_i\Delta\qquad i=1,\ldots,n-1;\cr
\Delta z_1&=q^{-1}z_n\Delta;\cr
\Delta H_{i{+}1}&=H_i\Delta\qquad i=1,\ldots,n-1;\cr
\Delta^2 H_1&=H_{n{-}1}\Delta^2.\cr}
$$
The last relation means that if we define $H_0:=\Delta
H_1\Delta^{-1}=\Delta^{-1}H_{n{-}1}\Delta$ then $H_0,\ldots,H_{n{-}1}$
generate an affine Hecke algebra while $H_1,\ldots,H_{n{-}1},\Delta$
generate the extended affine Hecke algebra corresponding to the weight
lattice of $GL_n$. In particular, there must be a family of $n$ commuting
elements: the Cherednik operators. In our particular case, they (or
rather their inverses) have a very nice explicit form. For $i=1,\ldots,n$
put
$$
\xi_i^{-1}:=\Hq_i\Hq_{i{+}1}\ldots\Hq_{n{-}1}\Delta H_1 H_2\ldots
H_{i{-}1}.
$$
We have the following commutation relations
$$
\eqalign{&\xi_{i{+}1}H_i=\Hq_i\xi_i;\quad
H_i\xi_{i{+}1}=\xi_i\Hq_i\qquad i=1,\ldots,n-1;\cr
&\xi_iH_j=H_j\xi_i\qquad i\ne j,j+1.\cr}
$$
The relation to Macdonald polynomials is as follows. For
$\lambda\in\Lambda$ define $\lamq\in k^n$ as
$\lamq_i:=q^{\lambda_i}t^{-k_i}$ where
$$
k_i:=\#\{j=1,\ldots,i-1\mid\lambda_j\ge\lambda_i\}
+\#\{j=i+1,\ldots,n\mid\lambda_j>\lambda_i\}.
$$
The following Lemma is easy to check:

\Lemma Self. {\rm(\cite{M1}~4.13,5.3)} a) The action of $\xi_i$ on $\cP$
is triangular. More precisely,
$\xi_i(z^\lambda)=
\lamq_iz^\lambda+\sum_{\mu<\lambda}c_{\lambda\mu}z^\mu$.\Par
b) Let $r\in\NN$. Then, with respect to the scalar product
$\<\cdot,\cdot\>_r$, the adjoints of $H_i$, $\Delta$, $\xi_i$ are
$H_i^{-1}$, $\Delta^{-1}$, $\xi_i^{-1}$ respectively.

Since $\lamq=\overline\mu$ implies $\lambda=\mu$, we immediately get:

\Corollary. The $\xi_i$ admit a simultaneous eigenbasis
$E_\lambda$ of the
form $z^\lambda+\sum_{\mu<\lambda}c_{\lambda\mu}z^\mu$ with
eigenvalue $\lamq_i$. Moreover, these functions coincide with those
defined in \cite{AAA}.

\noindent Actually, this gives a proof of \cite{AAA} and that the
$\xi_i$ commute pairwise.

\beginsection Recursion. Recursion relations

Observe that the operators $z_i$ and $\xi_i$ behave very similarly.
The only exception is, that there is no simple commutation rule for
$\Delta\xi_1$ as there is for $\Delta z_1$. Therefore, we introduce
another operator which simply is $\Phi:=z_n\Delta=q\Delta z_1$. Then
one easily checks the relations
$$
\Phi\xi_{i{+}1}=\xi_i\Phi;\quad \Phi\xi_1=q\xi_n\Phi.
$$
This implies:

\Theorem. Let $\lambda\in\Lambda$ with $\lambda_n\ne0$ and put
$\lambda^*:=(\lambda_n-1,\lambda_1,\ldots,\lambda_{n{-}1})$. Then
$E_\lambda=\Phi(E_{\lambda^*})$.

\noindent Observe that also the following relations hold:
$$
\eqalign{
\Phi z_{i{+}1}&=z_i\Phi\qquad i=1,\ldots,n-1;\cr
\Phi H_{i{+}1}&=H_i\Phi\qquad i=1,\ldots,n-1;\cr
\Phi^2 H_1&=H_{n{-}1}\Phi^2.\cr}
$$
This means that if $\tilde H_0:=\Phi
H_1\Phi^{-1}=\Phi^{-1}H_{n{-}1}\Phi$ then $\tilde
H_0,H_1,\ldots,H_{n{-}1}$ generate another copy of the affine Hecke
algebra, but note $\tilde H_0\ne H_0$! Indeed, $H_0$ is acting on
$\cP$, while $\tilde H_0$ acts only on $\cP'$.

The theorem above works as a recursion relation for $E_\lambda$ if
$\lambda_n\ne0$. The next (well known) lemma tells how to permute two
entries of $\lambda$.

\Theorem Haction. Let $\lambda\in\Lambda$ and $s_i$ a simple reflection.
\item{a)} Assume $\lambda_i=\lambda_{i+1}$. Then
$H_i(E_\lambda)=tE_\lambda$ and $\Hq_i(E_\lambda)=E_\lambda$.
\item{b)} Let $\lambda_i>\lambda_{i+1}$ and
$x:=1-\lamq_i/\lamq_{i{+}1}$. Then
$xE_\lambda=[xH_i+1-t]E_{s_i(\lambda)}$.

\Proof: a) Let $\mu<\lambda$. Then it follows from
properties of the Bruhat order that $\Hq_i(z^\mu)$ is a
linear combination of $z^\nu$ with $\nu<\lambda$. Hence, \cite{Self}a
implies that $H_i(E_\lambda)$ is orthogonal to all $z^\mu$ with
$\mu<\lambda$. By definition it must be proportional to $E_\lambda$.
The assertion follows by comparing the coeffients of $z^\lambda$.

b) Denote the right hand side by $E$. Since the coefficients of
$z^\lambda$ are the same, it suffices to prove that $E$ is an eigenvector
of $\xi_j$ with eigenvalue $\lamq_j$. This is only non\_trivial for
$j=i,i+1$. Let $j=i$. Then
$$
\eqalign{\xi_i(E)&=
\xi_i[x(\Hq_i+t-1)+(1-t)]E_{s_i(\lambda)}\cr
&=[x\xi_i\Hq_i+(1-t)\lamq_i/\lamq_{i{+}1}\xi_i]E_{s_i(\lambda)}=
[x\lamq_iH_i+(1-t)\lamq_i]E_{s_i(\lambda)}=\lamq_iE.\cr}
$$
The case $j=i+1$ is handled similarly.\qed

These two formulas suffice to generate all $E_\lambda$ but we will need
a more refined version.

\Lemma. Let $\lambda\in\Lambda$ with $\lambda_a\ne0$,
$\lambda_{a{+}1}=\ldots=\lambda_b=0$.
Let $\lambda^\#$ equal $\lambda$ except that the $a$-th and
$b$-th components are interchanged. Then
$$
(1-\lamq_at^a)E_\lambda=
[\Hq_a\Hq_{a{+}1}\ldots\Hq_{b{-}1}-\lamq_a t^a
H_aH_{a{+}1}\ldots H_{b{-}1}]E_{\lambda^\#}.
$$

\Proof: We prove this by induction on $b-a$. Let $\lambda'$ equal
$\lambda$ but with the $a$-th and $b-1$-st components interchanged. Then
$\lamq'_{b{-}1}=\lamq_a$ and $\lamq'_b=t^{-b+1}$. With
$x':=1-\lamq_a t^{b-1}$, \cite{Haction}b implies
$x'E_{\lambda'}=[x'H_{b{-}1}+1-t]E_{\lambda^\#}$. Hence,
with $x:=1-\lamq_at^a$, we get by induction
$$
\eqalign{
xx'E_\lambda&=
[\Hq_a\ldots\Hq_{b{-}2}-(1-x)H_a\ldots H_{b{-}2}]
[x'H_{b{-}1}+1-t]E_{\lambda^\#}\cr
&=[x'\Hq_a\ldots\Hq_{b{-}2}(\Hq_{b{-}1}+t-1)+
(1-t)\Hq_a\ldots\Hq_{b{-}2}-\cr
&-x'(1-x)H_a\ldots H_{b{-}1}-(1-t)(1-x)H_a\ldots H_{b{-}2}]
E_{\lambda^\#}\cr}
$$
Now observe $H_i(E_{\lambda^\#})=tE_{\lambda^\#}$ for
$i=a,\ldots,b-2$ and $\Hq_i(E_{\lambda^\#})=E_{\lambda^\#}$ by
\cite{Haction}a. Hence the expression above becomes
$$
\eqalign{[x'\Hq_a\ldots\Hq_{b{-}1}
-x'(1-x)&H_a\ldots H_{b{-}1}+\cr
&+(t-1)x'+(1-t)-(1-t)(1-x)t^{b-a-1}]
E_{\lambda^\#}.\cr}
$$
The lemma follows since the constant terms cancel out.\qed 

For $m=1,\ldots,n$ define the operators
$$
\eqalign{
A_m&:=H_mH_{m{+}1}\ldots H_{n{-}1}\Phi\cr
\Aq_m&:=\Hq_m\Hq_{m{+}1}\ldots \Hq_{n{-}1}\Phi\cr}
$$
Then we obtain:

\Corollary Recurs. Let $\lambda\in\Lambda$ with $m:=l(\lambda)>0$. Put
$\lambda^*:=
(\lambda_m-1,\lambda_1,\ldots,\lambda_{m{-}1},0,\ldots,0)$.
Then $(1-\lamq_mt^m)E_\lambda=[\Aq_m-\lamq_mt^m
A_m]E_{\lambda^*}$.

\beginsection Integrality. Integrality

To remove the denominators in the coefficients of $E_\lambda$ we use a
normalization as follows. Recall, that the {\it diagram} of
$\lambda\in\Lambda$ is the set of points (usually called {\it boxes})
$s=(i,j)\in\ZZ^2$ such that $1\le i\le n$ and $1\le j\le\lambda_i$. For
each box $s$ we define the {\it arm\_length} $a(s)$ and {\it
leg-length} $l(s)$ as
$$
\eqalign{
a(s)&\:=\lambda_i-j\cr
l'(s)&\:=\#\{k=1,\ldots, i-1\mid j\le \lambda_k+1\le\lambda_i\}\cr
l''(s)&\:=\#\{k=i+1,\ldots,n\mid j\le \lambda_k\le\lambda_i\}\cr
l(s)&\:=l'(s)+l''(s)\cr}
$$
If $\lambda\in\Lambda^+$ is a partition then $l'(s)=0$ and
$l''(s)=l(s)$ is just the usual leg\_length. Now we define
$$ 
\cE_\lambda:=\prod_{s\in\lambda}(1-q^{a(s)+1}t^{l(s)+1})E_\lambda.
$$
With this normalization, we obtain:

\Theorem Recurs2. With the notation of \cite{Recurs} let
$X_\lambda:=q^{\lambda_m-1}(\Aq_m-\lamq_mt^m A_m)$. Then
$\cE_\lambda=X_\lambda(\cE_{\lambda^*})$.

\Proof: It suffices to check the coefficient of $z^\lambda$. The factor
$q^{\lamq_m-1}$ cancels the effect of $\Phi=z_n\Delta$ on this
coefficient. The diagram of $\lambda^*$ is obtained from $\lambda$ by
taking the last non\_empty row, removing the first box $s_0$ and putting
the rest on top. It is easy to check that arm\_length and leg\_length of
the boxes $s\ne s_0$ don't change. Hence the assertion follows from
\cite{Recurs} since the factor corresponding to $s_0$ is just
$1-\lamq_mt^m$.\qed

Now we can state our first integrality result:

\Corollary Weak Integrality. Let $\cE_\lambda=\sum_\mu
c_{\lambda\mu}z^\mu$. Then $c_{\lambda\mu}\in\ZZ[t,q]$.

\Proof: Every $\cE_\lambda$ is obtained by repeated application of
operators $X_\mu$. Looking at the definition of $\Delta$ we conclude
that the $c_{\lambda\mu}$ are in $\ZZ[q,q^{-1},t]$. We exclude the
possibility of negative powers of $q$. For this write
$$
\Phi=z_n\Delta=z_n\Hq_{n{-}1}^{-1}\ldots\Hq_1^{-1}\xi_1^{-1}
$$
Now, $\cE_{\lambda^\#}$ is an eigenvector of $\xi_1^{-1}$ with
eigenvalue $(\lamq^*_m)^{-1}=q^{-\lambda_m+1}t^a$ with some
$a\in\ZZ$. This shows (by induction) that
$q^{\lamq_m-1}\Phi(\cE_{\lambda^*})$ doesn't contain negative powers
of $q$.\qed

Our goal is a more refined integrality result. For this we replace the
monomial basis $z^\lambda$ by a more suitable one. For
$\lambda\in\Lambda$ let $w_\lambda,\tilde w_\lambda$ be the shortest
permutations such that $\lambda^+:=w_\lambda^{-1}(\lambda)$ is
dominant (i.e. a partition) and $\lambda^-:=\tilde
w_\lambda^{-1}(\lambda)$ is antidominant. Now we define the {\it
$t$-monomial} $\fm_\lambda:=\Hq_{\tilde w^\lambda}(z^{\lambda^-})$.
The reason for this is that the action of the $H_i$ becomes nicer:
$$
H_i(\fm_\lambda)=\cases{
t\fm_{s_i(\lambda)}&if $\lambda_i\ge\lambda_{i{+}1}$\cr
\fm_{s_i(\lambda)}+(t-1)\fm_\lambda&if $\lambda_i<\lambda_{i{+}1}$\cr}
$$
$$
\Hq_i(\fm_\lambda)=\cases{
t\fm_{s_i(\lambda)}+(1-t)\fm_\lambda&if $\lambda_i>\lambda_{i{+}1}$\cr
\fm_{s_i(\lambda)}&if $\lambda_i\le\lambda_{i{+}1}$\cr}
$$
This is easily proved by induction on the length of $\tilde w_\lambda$.
Moreover, it is easy to see that the transition matrix between
$t$\_monomials and ordinary monomials is unitriangular.

Now we define a length function on $\Lambda$ by
$L(\lambda):=l(w_\lambda)=\#\{(i,j)\mid i<j,\lambda_i<\lambda_j\}$.

\Lemma. The function $\fm_\lambda^{(0)}=\sum_\mu t^{L(\mu)}\fm_\mu$,
where the sum runs through all permutations $\mu$ of $\lambda$, is
symmetric.

\Proof: It suffices to prove
$H_i(\fm_\lambda^{(0)})=t\fm_\lambda^{(0)}$ for all $i$. This follows
easily from the explicit description of the action given above.\qed

Clearly, the symmetric $t$\_monomials $\fm_\lambda^{(0)}$,
$\lambda\in\Lambda^+$ (later we will see that they are nothing else than
the Hall\_Littlewood polynomials) also have a unitriangular
transition matrix to the monomial symmetric functions $m_\lambda$.

For technical reasons we need also partially symmetric $t$\_monomials.
Let $0\le m\le n$ fixed and $\lambda\in\Lambda$ let
$\lambda':=(\lambda_1,\ldots,\lambda_m)$ and
$\lambda'':=(\lambda_{m{+}1},\ldots,\lambda_n)$. We also write
$\lambda=\lambda'\lambda''$. Let $\Lambda^{(m)}\subseteq\Lambda$ be
the set of those $\lambda$ such that $\lambda''$ is a partition. For
these elements we form
$$
\fm_\lambda^{(m)}:=\sum_\mu t^{L(\mu)}\fm_{\lambda'\mu}
$$
where $\mu$ runs through all permutations of $\lambda''$.

For $k\in\NN$ let $\phi_k(t):=(1-t)(1-t^2)\ldots(1-t^k)$. Then
$[k]!:=\phi_k(t)/(1-t)^k$ is the $t$\_factorial. For a partition $\mu$ we
define $m_i(\mu):=\#\{j\mid\mu_j=i\}$ and
$b_\mu(t):=\prod_{i\ge1}\phi_{m_i(\lambda)}(t)$. Now we define the
augmented partially symmetric $t$\_monomial as
$\tilde\fm_\lambda^{(m)}:=b_{\lambda''}(t)\fm_\lambda^{(m)}$.
The key result of this paper is

\Theorem Main1. For $m\ge l(\lambda)$ consider the expansion
$\cE_\lambda=
\sum\limits_{\mu\in\Lambda^{(m)}}c_{\lambda\mu}\tilde\fm_\mu^{(m)}$.
Then the coefficients $c_{\lambda\mu}$ are in $\ZZ[q,t]$.

\Proof: By \cite{Weak Integrality}, the only denominators which can
occur are products of factors of the form $1-t^k$ (or divisors thereof).
In particular, it suffices to show that the $c_{\lambda\mu}$ are in
$\ZZ[q,q^{-1},t,t^{-1}]$. Therefore, as in the proof \cite{Weak
Integrality}, we may replace $\Phi$ by 
$$
\Phi':=z_n\Hq_{n{-}1}^{-1}\ldots\Hq_1^{-1}=
t^{1-n}z_nH_{n{-}1}\ldots H_1
$$
and $A_m$, $\Aq_m$ by
$$
\eqalign{
A'_m&:=H_mH_{m{+}1}\ldots H_{n{-}1}\Phi'\cr
\Aq'_m&:=\Hq_m\Hq_{m{+}1}\ldots \Hq_{n{-}1}\Phi'\cr}
$$
Therefore, the theorem is proved with the next lemma.\qed

\Lemma OpInt. a) Every $\tilde\fm_\lambda^{(m)}$ is a linear combination
of $\tilde\fm_\mu^{(m+1)}$ with coefficients in $\ZZ[t]$.\Par
\noindent b) The operators $A_m'$ and $\Aq_m'$ commute with
$H_{m{+}1},\ldots,H_{n{-}1}$. In particular, they leave the space stable
which is spanned by all $\tilde\fm_\lambda^{(m)}$.\Par
\noindent c) The matrix coefficients of $A_m'$ and $\Aq_m'$ with respect
to this basis are in $\ZZ[t,t^{-1}]$.

\Proof: a) is obvious and b) an easy consequence of commutation and
braid relations. For c) let us start with another lemma.

\Lemma Phiop. Let $\lambda\in\Lambda$ with $\lambda_n\ne0$. Then
$\Phi'(\fm_{\lambda^*})=t^{-a}\fm_\lambda$ where
$a=\#\{i<n\mid\lambda_i>\lambda_n\}$.

\Proof: Using braid and commutation relations, one verifies
$\Phi'\Hq_i=\Hq_{i{-}1}\Phi'$ for all $i>1$. Now assume
$\lambda_{i{-}1}>\lambda_i$ for some $i<n$. Then using induction on
$l(\tilde w_\lambda)$ we may assume the result is correct for
$\fm_{s_{i{-}1}(\lambda)}$. Thus,
$$
\eqalign{
\Phi'(\fm_{\lambda^*})=&
\Phi'\Hq_i(\fm_{s_i(\lambda^*)})=
\Hq_{i{-}1}\Phi'(\fm_{s_i(\lambda^*)})=\cr
=&t^{-a}\Hq_{i{-}1}\fm_{s_{i{-}1}(\lambda)}=
t^{-a}\fm_\lambda\cr}
$$
Thus, we may assume $\lambda_1\le\ldots\le\lambda_{n{-}1}$. Let
$l:=\lambda_n-1$ and
$b:=\|max|\{i<n\mid\lambda_i\le l\}=n-1-a$. 
For simplicity, we denote $\fm_\lambda$ by $[\lambda]$. Thus
$$
\eqalign{
&t^{n{-}1}\Phi'[\lambda^*]=z_nH_{n{-}1}\ldots H_1[\lambda^*]
=t^b \Hq_{n{-}1}\ldots \Hq_{b{+}1}z_{b{+}1}
[\ldots,\lambda_b,l,\lambda_{b{+}1},\ldots]=\cr
&=t^b \Hq_{n{-}1}\ldots \Hq_{b{+}1}
[\ldots,\lambda_b,l+1,\lambda_{b{+}1},\ldots]
=t^b[\ldots,\lambda_{n{-}1},l+1]\cr}
$$\qed

\noindent We are continuing with the proof of \cite{OpInt}. By part a),
$\Aq'_m\fm_\lambda^{(m)}$ is symmetric in the variables
$z_{m{+}1},\ldots,z_n$. Therefore, it suffices to investigate the
coefficient of $[\mu]$ where $\mu''$ is an anti\_partition. Let $[\nu]$ be
a typical term of $\fm_\lambda^{(m)}$, i.e., $\nu'=\lambda'$ and $\nu''$
is a permutation of $\lambda''$. Let $l:=\lambda_1+1=\nu_1+1$. 
Looking how $\Phi'$ and $\Hq_i$ act, we see that $\Aq'_m[\nu]$ is a
linear combination of terms
$$
[\mu]=[s_{i_1}\ldots s_{i_r}(\nu_2,\ldots,\nu_n,l)]
$$
where $m\le i_1<\ldots<i_r<n$. If $r=n-m$ then
$b_{\mu''}=b_{\nu''}=b_{\lambda''}$ and we are done.

Otherwise there is $m\le j<n $ maximal such that
$j-1\not\in\{i_1,\ldots,i_r\}$. Since $s_j\ldots
s_{n{-}1}(\nu_2,\ldots,\nu_n,l)=(\ldots,\nu_j,l,\nu_{j{+}1},\ldots)$ we
necessarily have $\nu_j>l$. We are only interested in the case
where $\mu''$ is an anti\_partition. So $\nu_j$ has to be moved all the
way to the $m$\_th position. This means
$(i_1,\ldots,i_r)=(m,m+1,\ldots,j-2,j,\ldots,n-1)$. Moreover,
$\nu_m\le\ldots\nu_{j{-}1}\le l\le \nu_{j{+}1}\le\nu_n$ and $\nu_j>l$.
There are exactly $m_l(\nu'')+1=m_l(\mu'')$ such permutations $\nu''$
of $\lambda''$.
Each of them contributes $t^{-a}(1-t)t^{L(\nu'')}$ for the coefficient of
$[\mu]$ in $\Aq'_m[\lambda]$. With $j$, the
length $L(\nu'')$ runs through a consecutive segment
$b,\ldots,b+m_l(\nu'')$ of integers. So $[\mu]$ gets the factor
$t^{b-a}(1-t)(1+t+\ldots+t^{m_l(\nu'')})=
t^{b-a}(1-t^{m_l(\mu'')})$. This shows that the coefficient
of $[\mu]$ in $b_{\lambda''}A_m\fm_\lambda^{(m)}$ is divisible by
$b_{\mu''}$.

The case for $A_m$ is completely analogous and the
details are left to the reader. The only change is that one only
looks for the coefficient of $[\mu]$ where $\mu''$ is a partition.\qed

\noindent The proof gives actually a little bit more:

\Corollary q=0. For all $\lambda\in\Lambda$ we have
$\cE_\lambda(z;0,t)=\fm_\lambda$.

\Proof: Assume $m=l(\lambda)$ and look at $\Aq'_m\fm_{\lambda^*}$. In
the notation of the proof above, the second case (where some $s_j$ is
missing) can not occur since then $\nu_j=0$ would be greater than $l>0$.
So we get $\fm_\lambda=\Aq'_m\fm_{\lambda^*}$. This means that
$\fm_\lambda$ satisfies the same recursion relation as $\cE_\lambda$
with $q=0$.\qed

\beginsection Symmetric. The symmetric case and Kostka numbers

Finally we come to the integrality properties of the symmetric
polynomial $J_\lambda$ where $\lambda\in\Lambda^+$. For this we
normalize it as follows
$$
\cJ_\lambda(z;q,t)\:=
\prod_{s\in\lambda}(1-q^{a(s)}t^{l(s)+1})J_\lambda(z;q,t)
$$

\Theorem MainS. Let $\cJ_\lambda(z)=\sum_\mu
c_{\lambda\mu}(q,t)\tilde\fm^{(0)}_\mu$. Then the
coefficients $c_{\lambda\mu}$ are in $\ZZ[q,t]$.

\Proof: Let $m:=l(m)$ and consider $\lambda^-$, the anti\_partition with
$\lambda^-_i=\lambda_{n{+}1{-}i}$ and
$\lambda^0:=(\lambda_m-1,\ldots,\lambda_1-1,0,\ldots,0)$. Let
$\cE:=\Phi^m(\cE_{\lambda^0})$. Then $\cE$ equals $\cE_{\lambda^-}$,
except that in the normalization factor the contributions of the first
column of $\lambda^-$ are missing. Put
$$
\cJ:={(1-t)^m\over[m_0(\lambda)]!}\sum_{w\in W}H_w(\cE).
$$
We claim that $\cJ=\cJ_\lambda$. Consider the subspace $V$ of
$\cP$ spanned by all $E_{w\lambda}$, $w\in W$. Then it
follows from \cite{Self}b and the definitions that $\cJ_\lambda$ spans
$V^W$. This shows that $\cJ$ is proportional to $\cJ_\lambda$.

To show equality we compare the coefficient of $\fm_\lambda$. Since
$\lambda^-$ is anti\_dominant only those summands $H_w(\cE)$ have an
$\fm_\lambda$\_term where $w\lambda^-=\lambda$. These $w$ form a
left coset for $W_\lambda$. Therefore, summation over this coset
contribute the factor
$$
\sum_{w\in W_\lambda}t^{l(w)}=\prod_{i\ge0}[m_i(\lambda)]!=
{[m_0(\lambda)]!\over(1-t)^m}b_\lambda(t).
$$
Thus, $\fm_\lambda$ has in $\cJ$ the coefficient
$$
b_\lambda(t)\prod_{s\in\lambda^-\atop
s\ne(i,1)}\big(1-q^{a(s)+1}t^{l(s)+1}\big).
$$
On the other hand, by definition, the coefficient of $\fm_\lambda$ in
$\cJ_\lambda$ is
$$
\prod_{s\in\lambda}(1-q^{a(s)}t^{l(s)+1}).
$$
Let $w=w_{\lambda^-}$, the shortest permutation which transforms
$\lambda$ into $\lambda^-$. This means $w(i)>w(j)$ whenever
$\lambda_i>\lambda_j$ but $w(i)<w(j)$ for $\lambda_i=\lambda_j$ and
$i<j$. Consider the  following correspondence between boxes:
$$
\lambda\owns s=(i,j)\leftrightarrow s^-=(w(i),j+1)\in\lambda^- \ .
$$
This is defined for all $s$ with $j<\lambda_i$. One easily verifies
that $a(s)=a(s^-)+1$ and $l(s)=l(s^-)$. This means that $s$ and $s^-$
contribute the same factor in the products above. What is left out of
the correspondence are those boxes of $\lambda$ with $j=\lambda_i$ and
the first column of $\lambda^-$. The first type of these boxes
contributes $b_\lambda$ to the factor of $\cJ_\lambda$. The second
type doesn't 
contribute by construction. This shows $\cJ_\lambda=\cJ$.

Finally, we have to show that the coefficient of $\tilde\fm_\mu^{(0)}$ in
$\cJ$ is in $\ZZ[q,t]$. By \cite{Weak Integrality} it suffices to show
that these coefficients are in $A:=\ZZ[q,q^{-1},t,t^{-1}]$. So we can
ignore negative powers of $q$ and $t$ and replace $\cE$ by
$\cE':=(\Phi')^m(\cE_{\lambda^0})$ and similarly $\cJ$ by $\cJ'$.
\cite{Phiop} shows that $(\Phi')^m(\fm_{\mu}^{(m)})$ equals now the
symmetrization of some $\fm_{\mu'}$ in the {\it first\/} $n-m$ variables.
Therefore, by abuse of notation let $\mu=\mu''\mu'$ where $\mu'$ are the
last $m$ components of $\mu$ and let $\fm_\mu^{(m)}$ be the
symmetrization of $\fm_\mu$ in $\mu''$.

The isotropy group of $\lambda^-$ is generated by simple reflections.
Hence $\cE'$ is $W_{\lambda^-}$ invariant. We may assume
that $\mu$ is dominant for $W_\lambda$, i.e., $\mu_i\ge\mu_{i{+}1}$
whenever $\lambda_i=\lambda_{i{+}1}$. Let
$m_{ij}:=\#\{k\mid\lambda_k=i,\mu_k=j\}$. Then $\cE'$ is an
$A$\_linear combination of
$$
{b_{\mu''}\over \prod_{ij}[m_{ij}]!}\sum_{w\in
W_\lambda}H_w\fm_\mu $$
Averaging over $W$ we obtain that $\cJ'$ is an $A$\_linear combination
of
$$
{(1-t)^m\over[m_0(\lambda)]!}{b_{\mu''}\over \prod_{ij}[m_{ij}]!}
\prod_i[m_i(\lambda)]!\prod_j[m_j(\mu)]!\,\fm_\mu^{(0)}
$$
Because
$b_\mu=(1-t)^{n-m_0(\mu)}\prod_{j\ge1}[m_j(\mu)]!$ and
$b_{\mu''}=(1-t)^{l(\mu'')}\prod_{j\ge1}[m_{0j}]!$ the expression above
equals
$$
(1-t)^{m+l(\mu'')+m_0(\mu)-n}
{[m_0(\mu)]!\over[m_{00}]!}
\prod_{i\ge1}{[m_i(\lambda)]!\over\prod\limits_{j\ge0}[m_{ij}]!}
\,\tilde\fm_\mu^{(0)}.
$$
This proves the theorem.\qed

To put this into a more classical perspective note:

\Theorem. Let $\lambda\in\Lambda$. Then
$\fm_\lambda^{(0)}=\cJ_\lambda|_{q=0}$. Moreover, $\fm_\lambda^{(0)}$
(respectively $\tilde\fm_\lambda^{(0)}$) equals the Hall\_Littlewood
polynomial $P_\lambda(z;t)$ (respectively $Q_\lambda(z;t)$) in the
notation of \cite{M2}~III.

\Proof: The equality $\fm_\lambda^{(0)}=\cJ_\lambda\mid_{q=0}$
follows from \cite{q=0} and the equality $\cJ(z;0,t)=P_\lambda(z;t)$ is
well known (\cite{M2}~VI.1). Finally,
$\tilde\fm_\lambda^{(0)}=Q_\lambda(z;t)$ follows by comparing their
definitions.\qed

Recall that the Kostka\_functions $K_{\lambda\mu}(q,t)$ form the
transition matrix from the Macdonald polynomials $\cJ_\lambda(z;q,t)$ to
the $t$\_Schur functions $S_\mu(z;t)$. It is known that the transition
matrix from the $S_\mu(z;t)$ to the Hall\_Littlewood polynomials
$Q_\lambda(z;t)$ is unitriangular (\cite{M2}). Hence, \cite{MainS} can be
rephrased as

\Theorem. For all $\lambda,\mu\in\Lambda^+$, we have
$K_{\lambda\mu}(q,t)\in\ZZ[q,t]$.

\beginsection Jack. Jack polynomials

Let me shortly indicate how to obtain even positivity for Jack
polynomials. As already mentioned in the introduction, a detailed proof
completely in the framework of Jack polynomials will appear as a joint
paper with S.~Sahi. There we even give a combinatorial formula in terms
of certain tableaux.
\def\li#1{{\buildrel#1,\alpha\over\longrightarrow}}

Let $\alpha$ be an indeterminate and put formally $q=t^\alpha$. Let
$p(q,t)\in\QQ(q,t)$, $p_0\in\QQ$ and $k\in\NN$. Then we write $p\li{n}
p_0$ if $\|lim|\limits_{t\pfeil1}{p(t^\alpha,t)\over(1-t)^n}=p_0$.
For example, $1-q^at^b\li1a\alpha+b$.

Let $\cJ_\lambda(z;\alpha)$ be a Jack polynomial. One could define it by
$\cJ_\lambda(z;q,t)\li{|\lambda|}\cJ_\lambda(z;\alpha)$
(\cite{M2}~VI.10.23). There is also a non\_symmetric analogue defined
by $\cE_\lambda(z;q,t)\li{|\lambda|}\cE_\lambda(z;\alpha)$.

For any $\lambda\in\Lambda$ and $1\le m\le n$ we have
$\fm_\lambda\li0z^\lambda$ and
$\fm_\lambda^{(m)}\li0m_\lambda^{(m)}:=\sum_\mu z^{\lambda'\mu}$
where $\mu$ runs through all permutations of $\lambda''$.
If $\lambda\in\Lambda^+$ let $u_\lambda:=\prod_{i\ge1}m_i(\lambda)!$.
Then we have $b_\lambda\li{l(\lambda)}u_\lambda$. In particular,
$\tilde\fm_\lambda^{(m)}\li{l(\lambda)}\tilde
m_\lambda^{(m)}:=u_{\lambda''}m_\lambda^{(m)}$. With this notation we
have:

\Theorem. a) Let $\lambda\in\Lambda$ and $m=l(\lambda)$. Then there is
an expansion $\cE_\lambda(z;\alpha)=\sum_{\mu\in\Lambda^{(m)}}
c_{\lambda\mu}(\alpha)\tilde m_\mu^{(m)}$ with
$c_{\lambda\mu}(\alpha)\in\NN[\alpha]$ for all $\mu$.\Par\noindent
b) Let $\lambda\in\Lambda^+$ and $m=l(\lambda)$. Then there is
an expansion $\cJ_\lambda(z;\alpha)=\sum_{\mu\in\Lambda^+}
c'_{\lambda\mu}(\alpha)\tilde m_\mu^{(0)}$ with
$c'_{\lambda\mu}(\alpha)\in\NN[\alpha]$ for all $\mu$.

\Proof: Going to the limit $t\pfeil1$, \cite{Main1} and \cite{MainS}
imply $c_{\lambda\mu}, c'_{\lambda\mu}\in\ZZ[t]$ (see also
\cite{M2}~VI.10 Ex.~2a). It remains to show positivity. Since
$\cJ_\lambda(z;\alpha)$ is the symmetrization of
$\cE_\lambda(z;\alpha)$ it suffices to show positivity for the latter.

For this, we write the operator $X_\lambda$ of \cite{Main1} as
$X_\lambda=q^{\lamq_m-1}((\Aq_m-A_m)+(1-\lamq_mt^m)A_m)$. Then
$X_\lambda\li1X_\lambda^1$ and we prove that all parts of
$X_\lambda^1$ preserve positivity. With $\Phi_1:=z_ns_{n{-}1}\dots s_1$
this follows from
$$
\eqalign{
&q^{\lamq_m-1}\li01,\cr
&\Phi\li0\Phi_1,\cr
&A_m\li0s_m\ldots s_{n{-}1},\cr
&\Aq_m-A_m\li1\sum_{i=m}^{n-1}s_m\ldots\widehat{s_i}\ldots
s_{n-1}\Phi_1,\cr
&1-\lamq_mt^m\li1\alpha\lambda_m-k+m,\cr}
$$
where $k$ is number of $i=1,\ldots,m-1$ with $\lambda_i\ge\lambda_m$.
\qed

\beginsection References. References

\baselineskip12pt
\parskip2.5pt plus 1pt
\hyphenation{Hei-del-berg}
\def\L|Abk:#1|Sig:#2|Au:#3|Tit:#4|Zs:#5|Bd:#6|S:#7|J:#8||{%
\edef\TEST{[#2]}
\expandafter\ifx\csname#1\endcsname\TEST\relax\else
\immediate\write16{#1 hat sich geaendert!}\fi
\expandwrite\AUX{\neverexpand\ref{#1}{\TEST}}
\HI{[#2]}
\ifx-#3\relax\else{#3}: \fi
\ifx-#4\relax\else{#4}{\sfcode`.=3000.} \fi
\ifx-#5\relax\else{\it #5\/} \fi
\ifx-#6\relax\else{\bf #6} \fi
\ifx-#8\relax\else({#8})\fi
\ifx-#7\relax\else, {#7}\fi\Par}

\def\B|Abk:#1|Sig:#2|Au:#3|Tit:#4|Reihe:#5|Verlag:#6|Ort:#7|J:#8||{%
\edef\TEST{[#2]}
\expandafter\ifx\csname#1\endcsname\TEST\relax\else
\immediate\write16{#1 hat sich geaendert!}\fi
\expandwrite\AUX{\neverexpand\ref{#1}{\TEST}}
\HI{[#2]}
\ifx-#3\relax\else{#3}: \fi
\ifx-#4\relax\else{#4}{\sfcode`.=3000.} \fi
\ifx-#5\relax\else{(#5)} \fi
\ifx-#7\relax\else{#7:} \fi
\ifx-#6\relax\else{#6}\fi
\ifx-#8\relax\else{ #8}\fi\Par}

\def\Pr|Abk:#1|Sig:#2|Au:#3|Artikel:#4|Titel:#5|Hgr:#6|Reihe:{%
\edef\TEST{[#2]}
\expandafter\ifx\csname#1\endcsname\TEST\relax\else
\immediate\write16{#1 hat sich geaendert!}\fi
\expandwrite\AUX{\neverexpand\ref{#1}{\TEST}}
\HI{[#2]}
\ifx-#3\relax\else{#3}: \fi
\ifx-#4\relax\else{#4}{\sfcode`.=3000.} \fi
\ifx-#5\relax\else{In: \it #5}. \fi
\ifx-#6\relax\else{(#6)} \fi\PrII}
\def\PrII#1|Bd:#2|Verlag:#3|Ort:#4|S:#5|J:#6||{%
\ifx-#1\relax\else{#1} \fi
\ifx-#2\relax\else{\bf #2}, \fi
\ifx-#4\relax\else{#4:} \fi
\ifx-#3\relax\else{#3} \fi
\ifx-#6\relax\else{#6}\fi
\ifx-#5\relax\else{, #5}\fi\Par}
\setHI{[ABC]\ }
\sfcode`.=1000

\L|Abk:C1|Sig:C1|Au:Cherednik, I.|Tit:Double affine Hecke algebras
and Macdonald's conjectures|Zs:Annals Math.|Bd:141|S:191--216|J:1995||

\L|Abk:LS|Sig:LS|Au:Lascoux, A.; Sch\"utzenberger, M.%
|Tit:Sur une conjecture de H. O. Foulkes%
|Zs:C.R. Acad. Sci. Paris (S\'erie I)|Bd:286A|S:323--324|J:1978||

\L|Abk:M3|Sig:M1|Au:Macdonald, I.|Tit:A new class of
symmetric functions|Zs:Publ. I.R.M.A. Strasbourg, Actes
20$^{\hbox{e}}$ S\'eminaire Lotharingien|Bd:131-71|S:-|J:1988||

\B|Abk:M2|Sig:M2|Au:Macdonald, I.|Tit:Symmetric functions and Hall
polynomials (2nd ed.)|Reihe:-|Verlag:Clarendon Press|Ort:Oxford|J:1995||

\L|Abk:M1|Sig:M3|Au:Macdonald, I.|Tit:Affine Hecke algebras and
orthogonal polynomials|Zs:S\'eminaire Bourbaki|Bd:-|S:n$^o$ 797|J:1995||

\L|Abk:Op|Sig:O1|Au:Opdam, E.|Tit:Harmonic analysis for certain
representations of graded Hecke algebras%
|Zs:Acta Math.|Bd:175|S:75--121|J:1995||

\bye